\documentclass{elsart}

\newtheorem{Theoreme}{Thorme}[section]

\theoremstyle{definition}

\newtheorem{Remarque}[Theoreme]{Remark}

 \usepackage{graphicx}
 \usepackage{float}

\usepackage{amssymb}
\usepackage{amsmath}
\usepackage{placeins}

\def\rot{{\rm curl\,}}

\def\div{{\rm div\,}}

\def\Curl{{\rm Curl\,}}

\def\dblcaption#1#2{
\hbox to\textwidth{%
\hfill\vtop{\textwidth=12.0pc\hsize12.6pc\caption{#1}}\hfill\vtop{\textwidth=12.0pc%
\hsize12.6pc\caption{#2}}\hfill}}

\begin{document}

\begin{frontmatter}

\title{Simulation of the spread of a viscous fluid using a bidimensional shallow water model}
\author{Bernard Di Martino\corauthref{cor1}}
\ead{dimartin@univ-corse.fr}
\corauth[cor1]{Corresponding author:Tel. : +33 4 95 45 00 41}
\author{, Catherine Giacomoni,}
\author{Jean Martin Paoli}
\author{and Pierre Simonnet}%
\address{UMR CNRS 6134 SPE , Universit\'e de Corse, 20250 Corte, France}

\begin{abstract}
In this paper we propose a numerical method to solve the Cauchy
problem based on the viscous shallow water equations in an
horizontally moving domain. More precisely, we are interested in a
flooding and drying model, used to modelize the overflow of a
river or the intrusion of a tsunami on ground. We use a non
conservative form of the two-dimensional shallow water equations,
in eight velocity formulation and we build a numerical
approximation, based on the Arbitrary Lagrangian Eulerian
formulation, in order to compute the solution in the moving
domain.
\end{abstract}


\end{frontmatter}

\section{Introduction}

The flooding and drying of a fluid on ground is a problem which
has several applications in fluid mechanics such as coastal
engineering, artificial lake filling, river overflow or sea
intrusion on ground due to a tsunami. These works involves the
modelling of the physical process near the triple contact line
between the liquid (the water), the gas (the atmosphere) and the
solid (the ground). The numerical simulation of this problem is
very complex, due to the domain deformation and the triple contact
line movement.

Some numerical approaches to solve this problem have been
previously explored. For example, O. Bokhove \cite{bokhove1} uses
a conservative form of the shallow water equations associated to
the discontinuous galerkin discretization. D. Yuan {et al.}
\cite{yuan} propose to use a non conservative form of the shallow
water equations in total flow rate formulation and uses a finite
difference scheme for the numerical computation.

In this paper, we use a two dimensional viscous shallow water
model in eight velocity formulation. This model, called
Stokes-like by P.L. Lions (\cite{Lions}, section 8.3) is, in some
sense, intermediate between the semi-stationary model and the full
model of compressible isentropic Navier-Stokes model. This model
can be obtained by integrating vertically the so-called primitives
equations of the ocean \cite{LTW} with some hypothesis on the
viscous terms. We have chosen this model because we can prove the
existence of a solution under certain considerations as the
smoothness of the initial conditions and an acceptable hypothesis
on a boundary operator \cite{mateo}.

After the presentation of the model we describe the numerical
method based on the demonstration of the convergence proposed in
\cite{mateo}. Finally, we present some numerical examples in two
idealized domains and we give some perspectives for other studies.

\subsection{Mathematical model}

We assume that we know a continuous function $H$ from ${\mathbb
R}^2$ to ${\mathbb R}$ that represents the topography (the level of
the ground with respect to a reference level).

We note $\Omega_t$ (an open simply connected set of ${\mathbb
R}^2$) the horizontal domain occupied by the fluid at time $t$ and
$\eta(t,x,y)$ the elevation of the fluid compared to a horizontal
zero level. We set $h(t,x,y)$ the eight of the water column,
$h(t,x,y) = H(x,y) + \eta(t,x,y)$.

The shallow water equations are based on a depth integration of an incompressible fluid
conservation laws in a free surface-three dimensional domain. Governing equations for $u$ and $h$ can
be obtained in the usual way (\cite{berny1} for example) and the
two-dimensional system can be written as follows
\cite{mateo,mateo2,katia2}:
\begin{align}
&\frac{\partial u}{\partial t} + (u\cdot\nabla) u + g \nabla \eta
- \nu \Delta u + C_d |u| u = f_1,\qquad {\rm in } \bigcup_{t\in
[0,T]} \{ t\} \times \Omega_t  \label{eqqm}
\\ &\frac{\partial \eta}{\partial t} +
{\rm div} (uh)= 0,\qquad  {\rm in } \bigcup_{t\in [0,T]} \{ t\}
\times \Omega_t   \label{eqm}
\end{align}
where $g$ is the gravity constant, $\nu$ the viscosity
coefficient, $C_d$ a drag coefficient. $f_1$ represents
external forcing (for example the wind stress). These equations
are completed by initial conditions ($u(t=0)=u_0,\
\eta(t=0)=\eta_0$) and boundary conditions. In the following, we
assume that the boundary of the real domain is non-vertical. Hence
the variation of the surface elevation at the boundary imposes a
modification of the horizontal domain occupied by the fluid.

\subsection{Boundary conditions}

\subsubsection{In a three dimensional domain}

In many geophysical applications, the governing equations
are solved in a domain $[0,T]\times\Omega^3_t$ whereby
$\Omega^3_t=\{ (x,y,z) \in{\mathbb R^3} / (x,y)\in \Omega_t, -H(x,y)
< z < \eta(t,x,y)\}$ is a three dimensional moving domain. We denote by
$N=(n_1,n_2,n_3)$ the exterior unit normal to this domain, and we set
$N_h=(n_1,n_2)$ as its horizontal component. Moreover we denote
$V=(v_1,v_2,v_3)$ as the velocity of the fluid and $v_h=(v_1,v_2)$
as its horizontal component.

In the following, we assume that the ground of the sea does not have a
vertical wall. Then we have $n_3\neq 0$. The impermeability
boundary condition on the ground for the three dimensional problem
imposes $V(x,y,-H(x,y))\cdot N = 0$. Then, we have $v_h \cdot N_h = - v_3 n_3$ on the ground and, since
$n_3\neq 0$, we obtain $v_3 = -v_h \cdot \nabla H$.
\begin{figure}[H]
    \centerline{\includegraphics[width=11.truecm]{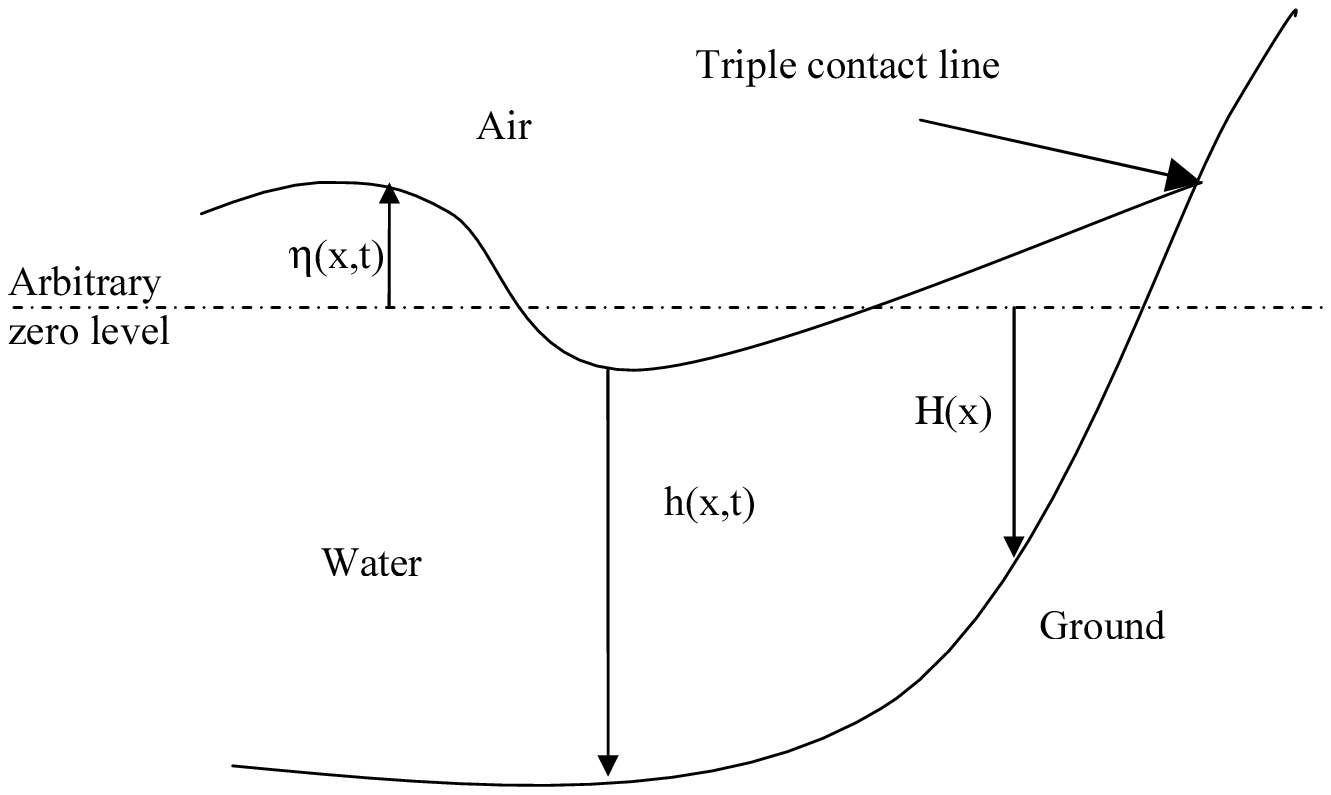}}
    \vskip -10truecm
    \caption{Vertical section of the three dimensional domain}\label{coupelaterale}
\end{figure}
At the surface, vertical velocity is equal to the lagrangian
derivative of the function $\eta$, $v_3 = \frac{\partial
\eta}{\partial t}+ v_h \cdot \nabla \eta$.

We recall that we are interested in the effect of the overflow
and we need to consider the problem in a moving domain that
represents the domain actually occupied by the fluid. Then we need
to impose a condition on the moving boundary. At the triple
contact line (between air, ground and water, see figure
\ref{coupelaterale}), the horizontal velocity at the surface
($v_h(x,y,\eta,t)$) is equal to the velocity at the bottom
($v_h(x,y,-H,t)$). Similarly, we have $v_3(x,y,\eta,t)=v_3(x,y,-H,t)$,
coherent with the fact that the lagrangian derivative of
$h$ is equal to zero on this boundary.

This property is not verified if, during its displacement, a part
of the triple contact line touches a vertical wall. In this case,
the boundary cannot be represented by the line $h=0$. This question is
very interesting but not treated in this paper.

\subsubsection{In a two dimensional domain}

The mean velocity $u$ used in the equations \ref{eqqm} and
\ref{eqm} is related to the previous three dimensional velocity by
$u=(u_1,u_2) = \frac{1}{h} \int_{-H}^\eta (v_1,v_2) {\rm d}z$
where $h\neq 0$. The two dimensional domain $\Omega_t$ is the
projection of the three dimensional domain $\Omega^3_t$ on the
horizontal plane.

For the depth integrated model, we denote by $\gamma_0$ the
boundary of $\Omega_0$ (at initial time). Assuming that
$\gamma_0$ is smooth enough, we define the deformed boundary as
follows : $\gamma_t=\{ x'=x+ d(x,t),\, x\in\gamma_0\}$, where
$d$ represents the horizontal displacement
$d(x,t)=\chi(0,t,x)-x$ and where $\chi(s,t,x)$ denotes the
Lagrangian flow, i.e. the position at time $t$ of the particle
located at $x$ at time $s$. With these notations, we have the
natural condition $\partial d(\chi(0,t,x),t)/\partial t=U(x,t)$ on
$\gamma_t$, where $U(x,t) = u(x+d(x,t),t)$.

%


The continuity equation (\ref{eqm}) can be written as follows:
\begin{equation}
h(\chi(x_0,0,t),t) = h(x_0,0) \exp \left(-\int_0^t \div u \right).
\label{eqdm1}
\end{equation}
So, if at the initial time, $h(x_0,0)=0$ then at all times
$h(\chi(x_0,0,t),t)=0$. Moreover, if $h(x_0,0)\neq 0$ then
$h(\chi(x_0,0,t),0)\neq 0 $ at all times. We conclude that the
boundary (and only the points of this boundary) at the initial
time can define the boundary at all times. For example, if at
initial time the domain is simply connected, the domain remains
simply connected at all times $t\geqslant 0$.

%

We characterize the boundary motion $\gamma_t$ by a condition on
the normal component of the fluid stress tensor $\sigma$ on the
boundary
\begin{equation} \label{tensionbord}
\sigma (x+d(x,t),t)\cdot n(x+d(x,t),t) | {\rm det} {J}| (x,t) =
A\left(\frac{\partial U(x,t)}{\partial t}\right),
\end{equation}
where $A$ is a differential operator defined on $\gamma_0$ taking
into account the boundary forcing on the fluid (see \cite{beal}
for example) and $J$ is the jacobian matrix associated to the
transformation $x\mapsto \chi(x,0,t)$. In \cite{mateo}, the
authors show that if the operator $A$ is a Laplace-Beltrami
operator ($\int_{\gamma_0} A(v)v = \int_{\gamma_0} A^{1/2}(v)
A^{1/2}(v) = |\!| v|\!|^2_{H^2(\gamma_0)})$ and if the data are
small enough, the problem has a solution. The condition on $A$ is
necessary to obtain the smoothness of the boundary and allows the
use of the classical Sobolev injections \cite{solo1,solo2}. More
precisely, in this paper, $D(A)=H^p(\gamma_0)$ with $p=2$ and
currently, we are not able to prove an existence result for $p<2$.
But numerically, it is possible to obtain physical acceptable
results with less smoothness on this operator (for example in
\cite{mateo2} we have obtained results by taking $A\equiv 0$).

It is not easy to characterize the operator $A$. It needs a good
physical interpretation and give sufficient mathematical
smoothness to prove the existence of a solution of our problem.
Physically, this operator represents a friction condition in the
triple contact line (between the fluid, the solid domain and the
air). In the literature, a large number of physical approaches in
the case of a wetting film of fluid or for a droplet
\cite{DeGennes} are found. Unfortunately, we cannot use directly
these results. Firstly, we do not work at the microscopic scale
(mainly governed by the Van Der Waals forces) and we need to take
into account the effect of these phenomena at large scale.
Secondly, the main effect studied for those phenomenon is in the
vertical direction (for example for the capillary effect) but
here, we use a depth integrated model.

\section{The numerical method}
\subsection{The ALE equations}

We use a numerical method based on the weak formulation of the
problem in the domain $\Omega_t$ depending on time. The ALE method that here we use, has been used to solve
 the Navier-Stokes equations in a moving domain \cite{maury}.

We assume that at the initial time $t_0$, the fluid domain is
covered by a regular mesh (for example, an usual finite element
mesh). We also assume that on the boundary, each point of this
mesh has the velocity of the fluid. The velocity of the internal
points can be chosen arbitrarily, the only condition is to
conserve the smoothness of the mesh at all times.

Then, at any time $t \in [0,T]$, the velocity $c^t$ of the moving
mesh can be defined by solving the following problem:
\begin{align}
& \Delta c^t = 0 \mbox{ in } \Omega_t \label{pbC1} \\
& c^t =u(x,t)\ \text{on}\ \gamma_t \label{pbC2}
\end{align}
and we consider the following mapping
\begin{align*}
c: \tilde Q &\rightarrow {\mathbb R}^2 \\
 (x,t) &\mapsto c(x,t) = c^t(x),
\end{align*}
where $\displaystyle \widetilde{Q}=\cup_{t\in
[0,T]}(\Omega_t\times\{t\})$. The relation between the different
domains $\Omega_{t_i}$ is given by the mapping
\begin{align*}
{\mathcal C}(\cdot, t_1,t_2): \Omega_{t_1} &\rightarrow
\Omega_{t_2}
\\
 x_1&\mapsto x_2={\mathcal C}(x_1,t_1,t_2)
\end{align*}
where $\mathcal C(.,t_1,t_2)$ is the characteristic curve from
$(x_1,t_1)$ to $(x_2,t_2)$ corresponding to the velocity $c^t$ in
the space-time domain $\displaystyle \widetilde{Q}$. For each time
$\tau$, we denote by $u_\tau$ and $h_\tau$ the ALE velocity and
the ALE thickness:
\begin{equation}\label{eq1}
u_\tau(x,t) = u({\mathcal C}(x,\tau,t),t), \mbox{ } h_\tau(x,t) =
h({\mathcal C}(x,\tau,t),t).
\end{equation}
With these notations, we can write the ALE formulation of the
problem
\begin{align} & \frac{\partial u_\tau}{\partial t} +
(u_\tau - c_\tau)\nabla u_\tau - \mu \Delta u_\tau + C_d\  u_\tau
| u_\tau |+ g\nabla h_\tau = O(t-\tau),
\\
& \frac{\partial h_\tau}{\partial t} + (u_\tau - c_\tau)\nabla
h_\tau + h_\tau\div u_\tau = O(t-\tau).
\end{align}

In what follows, we will use a first order time discretization to
solve this problem (see \cite{Fourestey} for more details about
second order schemes).

\subsection{Time discretization}

Let $M\geqslant 1$, we note $\Delta t=\frac{T}{M}$ as the time step.
For $1\leqslant n\leqslant M$, $1\leqslant m\leqslant M $, we set
$\Omega^n=\Omega_{t^n}$ with boundary $\gamma^n$ and
\begin{align*}
g_m^n(x) &= g_{t^m}(x,t^n)\ \text{for}\ x\in \Omega^m, \\ g^n(x) &=
g_n^n(x)=g(x,t^n)\ \text{for}\ x\in \Omega^n.
\end{align*}
In order to ensure the positivity of $h_{t^n}$, the continuity
equation is renormalised as follows
\begin{equation}\frac{\partial}{\partial t} \log h_{t^n} +
(u_{t^n}-c_{t^n})\nabla \log h_{t^n} + \div u_{t^n} = 0.
\end{equation}
With this renormalisation, we do not have conservation of the
mass, but when we follow the evolution of the mass during a
simulation (e.g. in the first test presented in the following
section), only small variations of the mass quantity are observed.

We denote by $U^n$ (respectively $H^n$ and $C^n$) the
approximation of the exact solution $u_{t^n}$ (respectively
$h_{t^n}$ and $c_{t^n}$).

\noindent Then, we note $\tilde U^n(x) = U^n(X^n(x,t^n))$ and
$\tilde H^n(x) = H^n(X^n(x,t^n))$ with $X^n(x,\cdot)$ the
characteristic curve solution of: $$ \left\{
\begin{array}{l}
\displaystyle \frac{\partial X^n}{\partial t}(x,t) =
(U^n-C^n)(X^n(x,t))\\ X^n(x,t^{n+1}) = x.
\end{array}
\right. $$ In our study, we approximate the foot of the
characteristic by $X^n(x,t^n) \simeq x - (U^n-C^n)(x) \Delta t.$
We set $U_n^{n+1}$ (resp. $H_n^{n+1}$) the approximation of
$u_{t^n}(x,t^{n+1})$ (resp. $h_{t^n}(x,t^{n+1})$. Approximating
the Lagrangian derivative by a first order Euler scheme, and using
a linearised drag operator, we obtain the following:
\begin{align}
&U_n^{n+1}+\Delta t\, \left( g \nabla H_n^{n+1} +  C_d\ U_n^{n+1}
| \tilde U^n| - \nu\, \Delta U_n^{n+1} \right) = \tilde U^n,&
\text{in}\ \Omega^n \label{eq11}
 \\
& \log H^{n+1}_n + \Delta t\, \div U_n^{n+1} = \log \tilde H^n &
\text{in}\ \Omega^n,\label{eq12} \\ &\div U_n^{n+1} =
A\left(\frac{\partial{\tilde U^n}}{\partial t}\right) & \text{on}\
\gamma^n,\\ &\rot U_n^{n+1} =  H_n^{n+1} = 0& \text{on}\ \gamma^n.
\end{align}
Implicitly, taking into account the equation \ref{eqdm1}, we have
$H_n^{n+1}=0\ \text{on}\ \gamma^n$. This is an approximation at
the first order in time of the initial shallow water problem.

\begin{Remarque}
\item This problem is solved by a fixed point technic. A first
approximation of \ref{eq12} is computed assuming $U_n^{n+1} =
\tilde U^{n}$ and the approximation of $H^{n+1}_n$ is then used in
\ref{eq11} to obtain a first approximation of $U_n^{n+1}$. We
repeat this operation in order to obtain convergence of the
system. Prove of convergence can be fount in \cite{mateo}.

\item We cannot easily increase the order of this scheme because
the ALE formulation is only of the first order.
\end{Remarque}

The previous problem is solved on $\Omega^n$ by using the spatial
discretization proposed in the following section. After, we need to solve the problem given by equations
\ref{pbC1} and \ref{pbC2} in order to compute the mesh velocity.
Each point $p_k$ with coordinates $x_k$ of the mesh is then moved
with the first order approximation $x_k^{n+1} = x_k^n + \Delta t
C^n(x_k^n,t^n)$.

\subsection{Spatial discretization}

The spatial discretization of the previous problem is based on the
Finite Element Galerkin method. In  \cite{mateo}, an approach
based on the Galerkin method with a special basis is proposed, but
this approach cannot be applied easily if $A \not\equiv 0$.

We note $\{\phi_k\}$ the set of finite element functions of $H^1(\Omega^n)$. The weak
formulation of our approximated problem is
\begin{align}
&\int_{\Omega^n}U_n^{n+1}(x) \phi_k(x) + \nu\,\Delta t\,
\int_{\Omega^n} \nabla U_n^{n+1}(x) \nabla  \phi_k(x) \notag\\ & -
g\,\Delta t\, \int_{\Omega^n} H_n^{n+1}(x) \div \phi_k(x) +
C_d\,\Delta t\, \int_{\Omega^n} U_n^{n+1}(x) | U_n^{n}(x)|
\phi_k(x) \notag
\\ & + \int_{\gamma_t} \sigma (X+d(X,t),t)\cdot n(X+d(X,t),t)
\gamma (\phi_k(x))  \notag \\ & = \int_{\Omega^n} \tilde U^n(x)
\phi_k(x) & \label{equa_10}
 \\
& \log H^{n+1}_n(x) + \Delta t\, \div U_n^{n+1}(x) = \log \tilde
H^n(x) &  \text{in}\ \Omega^n
\end{align}
where the operator $\gamma$ is the trace operator from $\Omega^n$
to $\gamma^n$.

Then the equations (\ref{tensionbord}) and (\ref{equa_10}) give
\begin{align}
&\int_{\Omega^n}U_n^{n+1}(x) \phi_k(x) + \nu\,\Delta t\,
\int_{\Omega^n} \nabla U_n^{n+1}(x) \nabla  \phi_k(x) \notag\\ & -
g\,\Delta t\, \int_{\Omega^n} H_n^{n+1}(x) \div \phi_k(x) + C_d
\int_{\Omega^n} U_n^{n+1}(x) | U_n^{n}(x)| \phi_k(x) \notag
\\ & + \int_{\gamma_0} A^{1/2}\left(\frac{\partial \tilde
U}{\partial t}\right) A^{1/2} \left(\gamma (\tilde \phi_k(x_0))\right) d
\gamma_0 \notag
\\ & = \int_{\Omega^n} \tilde U^n(x) \phi_k(x)
\end{align}
where $\tilde g (x) = g(\chi(x,0,t))$.

We then use a first order discretization of the partial time
derivative:
\begin{align}
&\int_{\gamma_0} A^{1/2}\left(\frac{\partial \tilde U}{\partial
t}\right) A^{1/2} \left( \gamma (\tilde \phi_k(x_0)) \right) d
\gamma_0 \simeq \frac{1}{\Delta t}\int_{\gamma_0} A^{1/2}\left(
\tilde U_n^{n+1}(x)\right) A^{1/2} \left(\gamma (\tilde
\phi_k(x_0))\right) d \gamma_0 \notag\\ & - \frac{1}{\Delta
t}\int_{\gamma_0} A^{1/2} \left(\tilde U_n (x)\right)
 A^{1/2} \left(\gamma (\tilde \phi_k(x_0))\right) d
\gamma_0.
\end{align}

The first part of the right hand of this equation needs to be
included on the left hand of the global problem and the second
part on the right hand.

The global weak formulation is:
\begin{align}
&\int_{\Omega^n}U_n^{n+1}(x) \phi_k(x) + \nu\,\Delta t\,
\int_{\Omega^n} \nabla U_n^{n+1}(x) \nabla  \phi_k(x)  - g\,\Delta
t\, \int_{\Omega^n} H_n^{n+1}(x) \div \phi_k(x) \notag\\ &
\int_{\Omega^n} U_n^{n+1}(x) | U_n^{n}(x)| \phi_k(x) +
 \frac{1}{\Delta t}\int_{\gamma_0} A^{1/2} \left(\tilde U_n^{n+1}(x)\right)
A^{1/2} \left(\gamma (\tilde \phi_k(x_0))\right) d \gamma_0 \notag
\\ & = \int_{\Omega^n} \tilde U^n(x) \phi_k(x) +  \frac{1}{\Delta t}\int_{\gamma_0}
A^{1/2}\left(\tilde U_n (x)\right) A^{1/2} \left(\gamma (\tilde \phi_k(x_0))\right) d
\gamma_0, & \label{equa_11}
 \\
& \log H^{n+1}_n(x) + \Delta t\, \div U_n^{n+1}(x) = \log \tilde
H^n(x). &
\end{align}

\begin{Remarque}
On $\gamma_t$, $U(l,t) = \sum_k a_k \phi_k(l,t)$ where $l$ is a
point of $\gamma_t$ and the sum is computed on all finite element
functions. Then, on the initial boundary $\gamma_0$, we have the
decomposition $\tilde U(l_0,t) =\sum_k a_k \tilde \phi_k(l_0,t)$
with the same $a_k$.
\end{Remarque}

\begin{Remarque}
Formally, we can use two kinds of boundary conditions:
\begin{enumerate}
\item A condition of normal displacement of the boundary. We
write $u\cdot n = \frac{\partial d}{\partial t}$ where $d$ is the
normal displacement. With this condition, the well posed weak
formulation associated to the diffusion operator is $\int_\Omega
\div u \div \phi + \int_\Omega \rot u \rot \phi - \int_\gamma \div
u \phi\cdot n + \int_\gamma \rot u \phi\cdot\alpha(n)$ where
$\alpha(u_1,u_2) = (-u_2,u_1)$ and $\rot u = \frac{\partial
u_2}{\partial y} - \frac{\partial u_1}{\partial x}$. The boundary
condition is then imposed on $h - \div u$ and $\rot u$ since we
do not want to impose a condition on $u\cdot n$ or $u\cdot \alpha(n)$.
\item A condition on the displacement of the boundary  $u = \frac{\partial d}{\partial
t}$ where $d$ is a vector. With this condition, the well posed
weak formulation associated to the diffusion operator is
$\int_\Omega \nabla u \nabla \phi  - \int_\gamma \nabla u
\phi\cdot n$ and the boundary condition is then imposed on $\nabla
u$ (more precisely, we need to take into account the trace
tensor including condition on $h$ in the boundary, but we assume
$h\equiv0$ on this boundary).
\end{enumerate}
\end{Remarque}
\begin{Remarque}  \label{remarqueoscil}
 We can note that the friction term $C_d u |u|$ is necessary to
stabilize the flow. Indeed, if we use the following decomposition
\begin{equation}
(L^2(\Omega^n))^2=\nabla H_0^1(\Omega^n)\oplus \Curl
H_0^1(\Omega^n)\oplus\nabla {\mathcal H}(\Omega^n),
\end{equation}
where ${\mathcal H}$ is the intersection of $H^1$ and the space of
harmonic functions, we can see that the functions of $\nabla
{\mathcal H}(\Omega^n)$ are not controlled by the laplacian
operator. We present in the following section a case with $C_d=0$
where we do not have damping of the flow.
\end{Remarque}

\subsection{First numerical test}

\begin{figure}[ht]
    \centerline{\includegraphics[width=11.truecm, angle=-90]{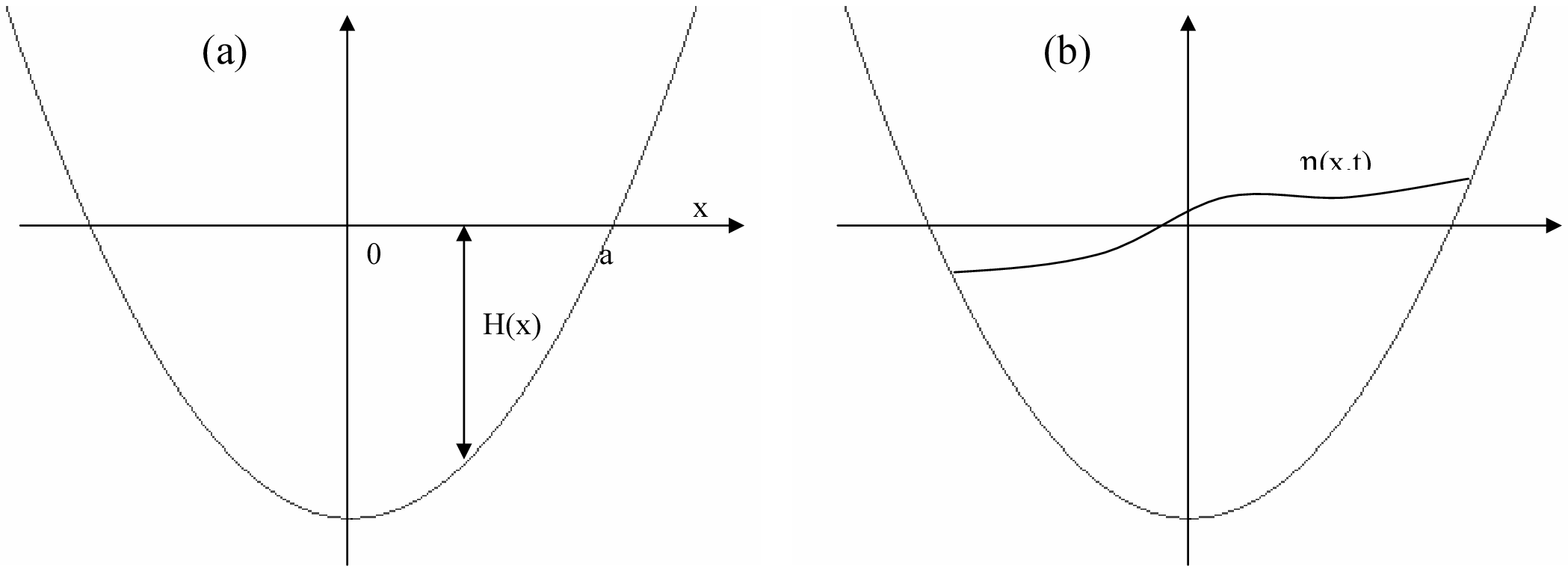}}
    \vskip -5truecm
    \caption{Simplified domain - (a) at rest. (b) at an arbitrarily time.}\label{domainvert}
\end{figure}

The previous method is tested in order to simulate the behavior of
a fluid in a simplified domain. We assume that the domain is
axisymmetric using the function $H(r) :
r\mapsto a r^2+1$. At the initial time $t=0$, the fluid at rest
touches the wall ($h(a,t=0)=0$) (see figure \ref{domainvert}).
Physically, our initial domain $\Omega_0$ is a disc with a radius
of 130 meters. This domain is meshed with triangles (420 triangles
for the mesh M1, 470 for the mesh M2 and 1002 for the mesh M3).
The height of the column of water at the centre is 1 meter (since $r=0$
at the centre of our domain). Fluid viscosity is $0.01 m^2.s^{-1}$
and gravity coefficient $g$ is assumed equal to $1 m.s^{-2}$ in
order to amplify the elevation.

 We apply a forcing at the surface of this fluid. This
forcing is usually taken into account by applying the continuity
of the horizontal stress tensor on the surface. So, for the three
dimensional model, $$\frac{\partial u_h}{\partial z} = C W_a |W_a|
$$ where $W_a$ is the wind velocity (at ten meters over the flow
for the ocean for example) and $C$ is a ``drag coefficient''
assumed to be constant. Using the vertical integration of the
vertical operator $\frac{\partial}{\partial z} \nu_z
\frac{\partial u_h}{\partial z}$, and assuming $u_h(\eta)\simeq
u$, we obtain a forcing condition on the mean velocity $f \propto
\frac{ W_a |W_a|}{H}$.

\begin{figure}[ht]
    \centerline{\includegraphics[width=15.truecm]{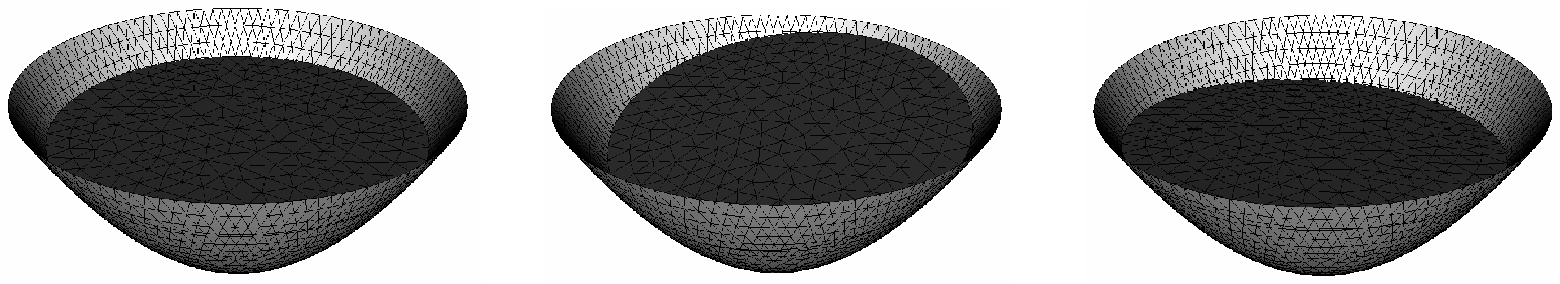}}
    \vskip -18truecm
    \caption{Oscillation of the water}\label{vuecoupe}
\end{figure}

\begin{figure}[ht]
    \centerline{\includegraphics[width=7.truecm, angle=-90]{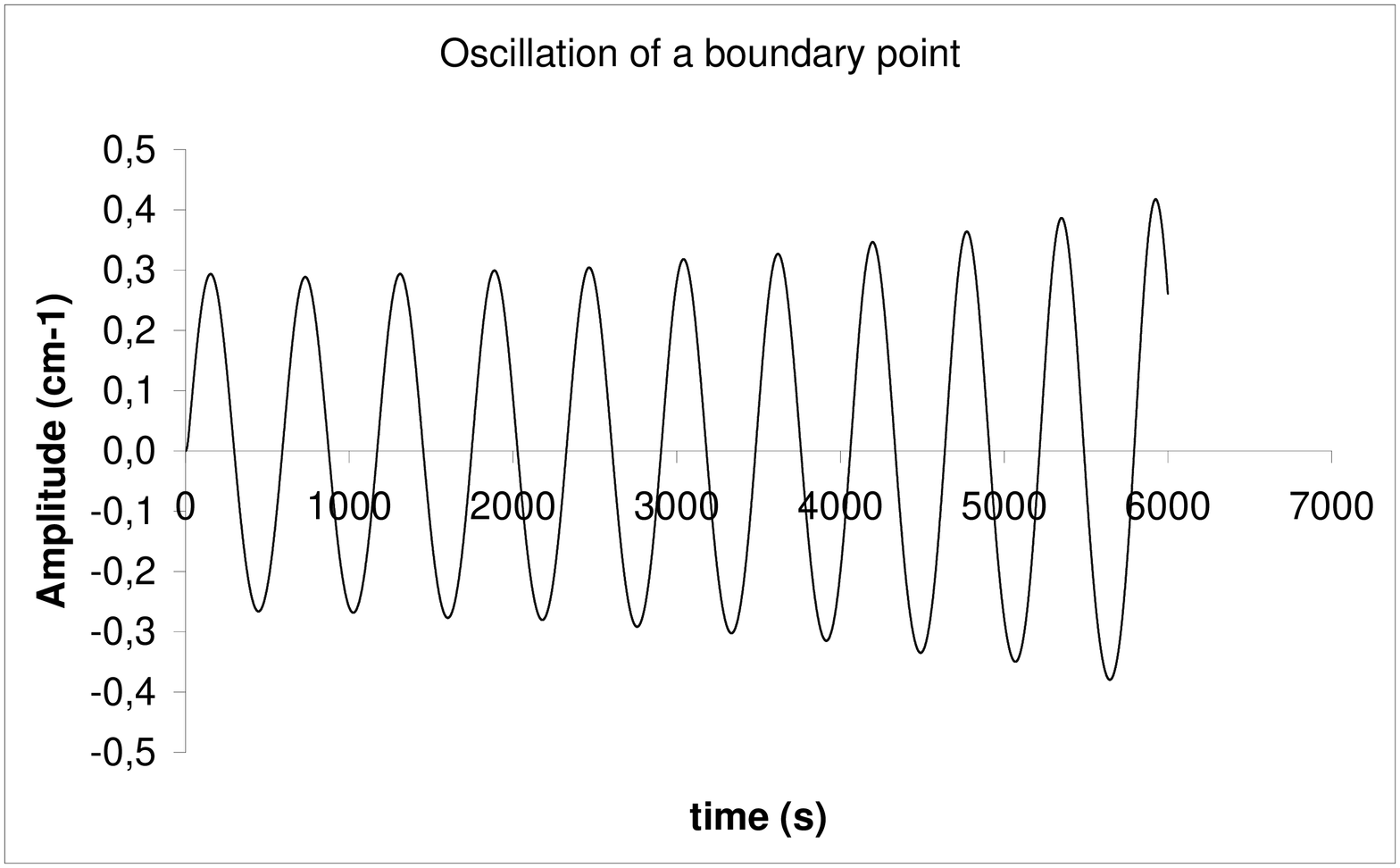}}
    \caption{Oscillations of the water without drag term}\label{oscillations}
\end{figure}

We use here
\begin{equation} \label{forcing}
\left\{
\begin{array}{cc}
 f = 1 &\ if\ 0<t<20\ {\rm s} \\
 f = 0 &\ if \ t>20\ {\rm s} \\
\end{array}.\right.
\end{equation}
Hence, we observe an oscillation of the free surface of the water
plan. At each oscillation, a part of the water moves on ground and
a part of the ground is uncovered by this water (see figure
\ref{vuecoupe}).
In the last part of the remark \ref{remarqueoscil}, we indicate
that with $C_d=0$ a part of the flow component is not diffusive.
More precisely, our solution is only a gradient. If we do not take
into account boundary friction effects or bottom friction effects,
a part of the fluid is not diffusive. To observe this effect, we
plot (on figure \ref{oscillations}) the level of a boundary point
according to time.
We can observe an accumulation of energy on the fluid due to the
numerical approximation, and this accumulation is not compensated
by diffusive term.

\begin{figure}[ht]
    \centerline{\includegraphics[width=7.truecm, angle=-90]{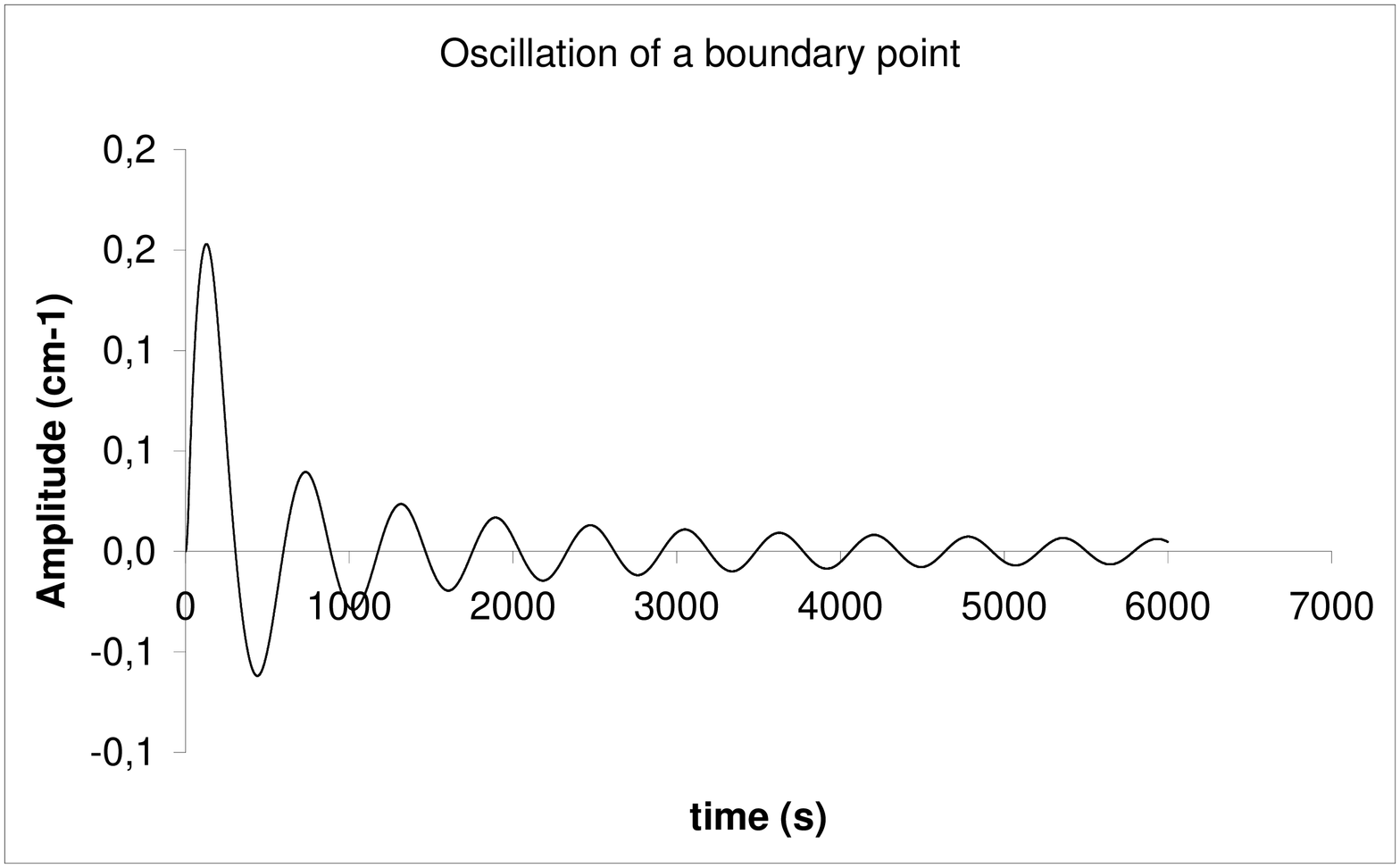}}
    \vskip -0.5truecm
    \caption{Oscillations of the water with a drag term}\label{oscillationsamorties}
\end{figure}

In figure \ref{oscillationsamorties}, we plot the variation of the
same point taking into account the drag coefficient $C_d |u| u$.
With this term, all the modes of the fluid are diffusive and, if we
stop the forcing, global energy of the fluid vanishes.

\begin{figure}[ht]
    \centerline{\includegraphics[width=7.truecm, angle=-90]{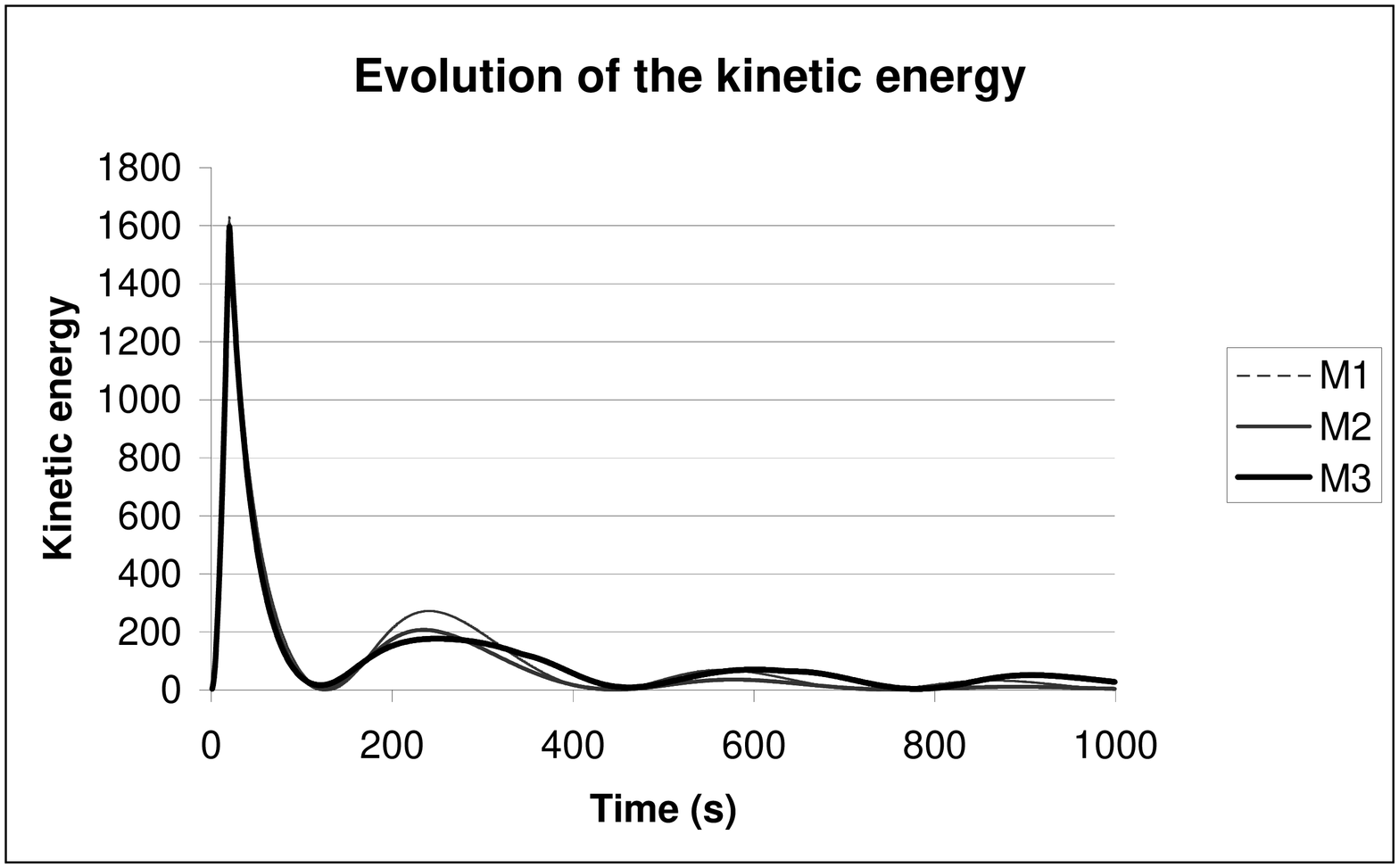}}
     \vskip -0.5truecm
    \caption{Evolution of the solution for M1, M2 and M3 meshes}\label{evolutionkinetic}
\end{figure}

In figure \ref{evolutionkinetic}, we plot the time evolution of the
kinetic energy of the flow for the different meshes (M1, M2 and M3).
We only have a little difference for all these meshes but a very
expensive cost for the mesh M3.
After the forcing phase (20 seconds), we can observe the
transformation of the kinetic energy in potential energy and
conversely. When the kinetic energy vanishes, potential energy is
maximal. If the drag coefficient is assumed to be equal to zero,
the amplitude of the oscillation does not decrease.

\begin{figure}[ht]
    \centerline{\includegraphics[width=7.truecm, angle=-90]{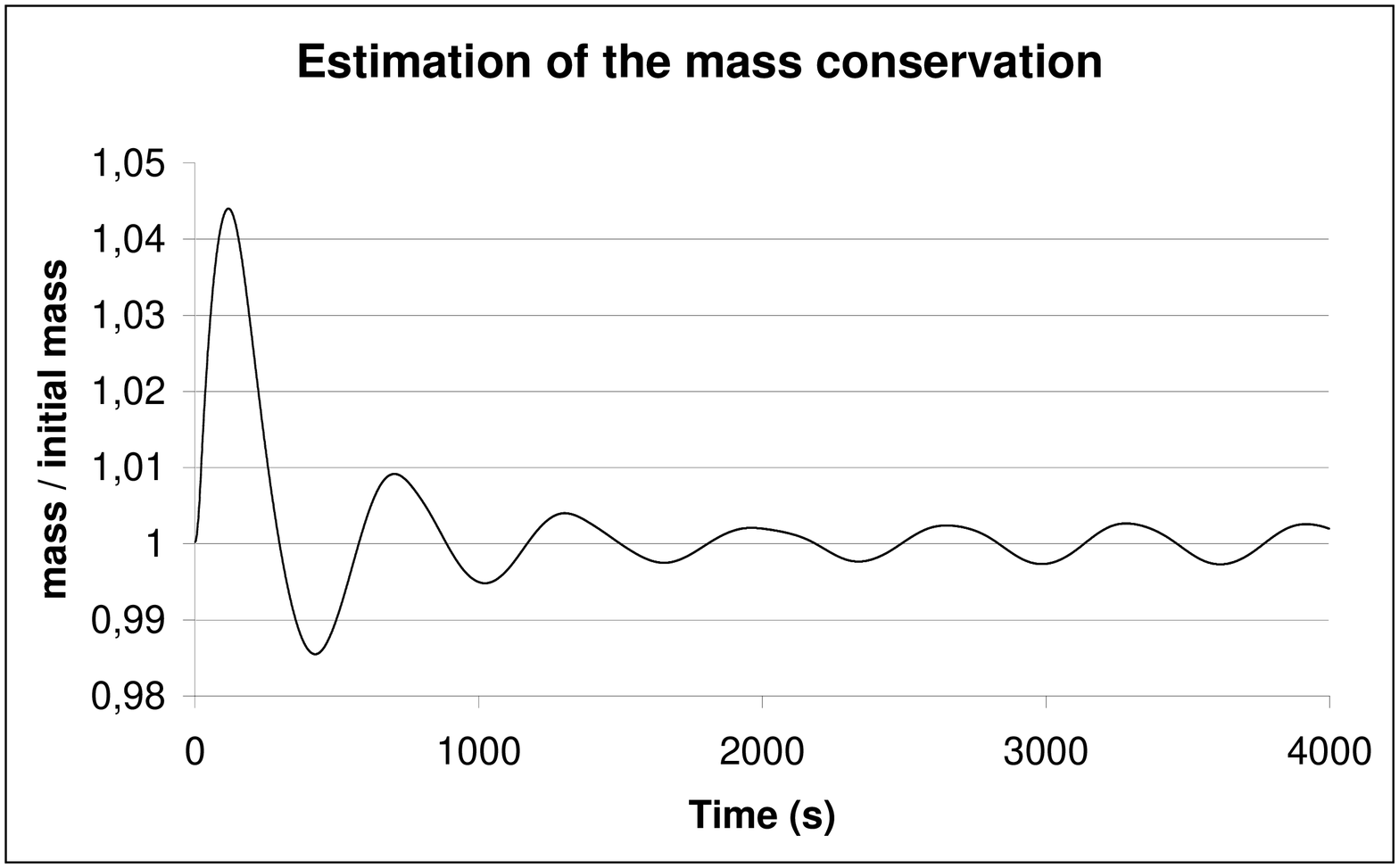}}
    \vskip -0.5truecm
    \caption{Evolution of the mass in comparison with initial mass}\label{massconservation}
\end{figure}

Finally, figure \ref{massconservation} represents the evolution of
the mass of fluid compared to the initial mass. Even if we do not
have a rigorous conservation of the mass, due to the
renormalisation of the mass equation, the variation of the mass is
very little.

\FloatBarrier

\subsection{Second numerical test}

\begin{figure}[ht]
    \centerline{\includegraphics[width=9.truecm, angle=-90]{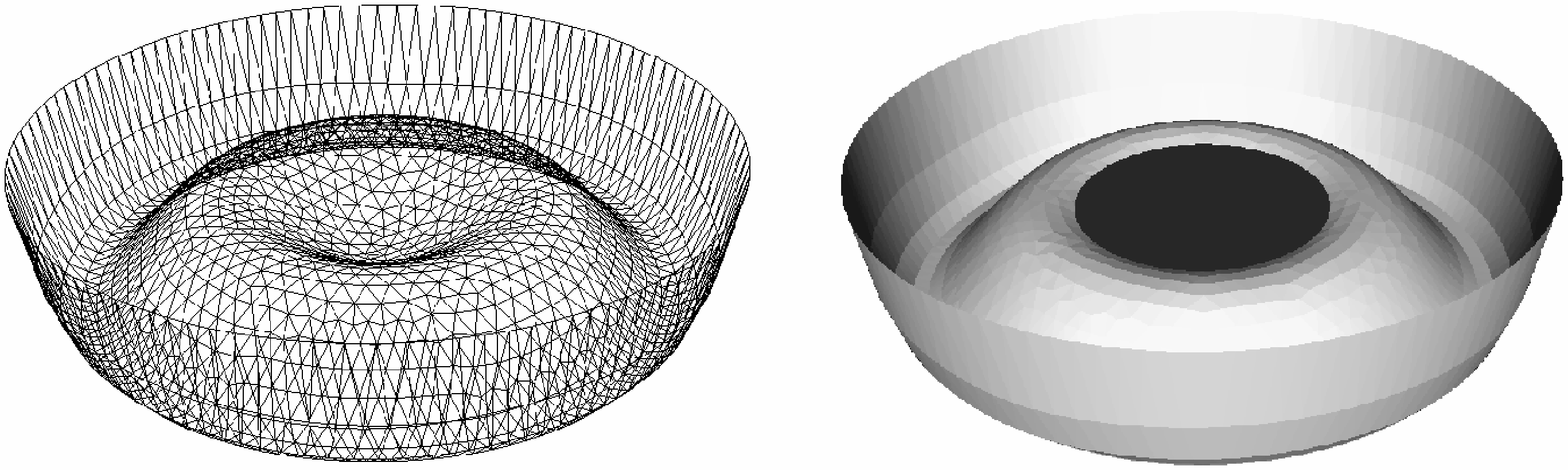}}
    \vskip -2truecm
    \caption{Mesh of the ground and initial position of the water}\label{Dom_ini}
\end{figure}

In this second numerical experiment, we use a more complex domain
to test our numerical method. We assume that the domain is
axisymmetric using the function $H(r) :
r\mapsto a r^3+ b r^2 + c r+d$. $b,c$ and $d$ are chosen in order
to at the distance $1/a$ of the centre of the fluid domain, the
fluid touches the wall (see figure \ref{Dom_ini}).
We use the same external forcing as for the first experiment
(\ref{forcing}). But, due to the specific form of the topography,
a part of the fluid goes to the external crown.

\begin{figure}[ht]
    \centerline{\includegraphics[width=9.truecm, angle=-90]{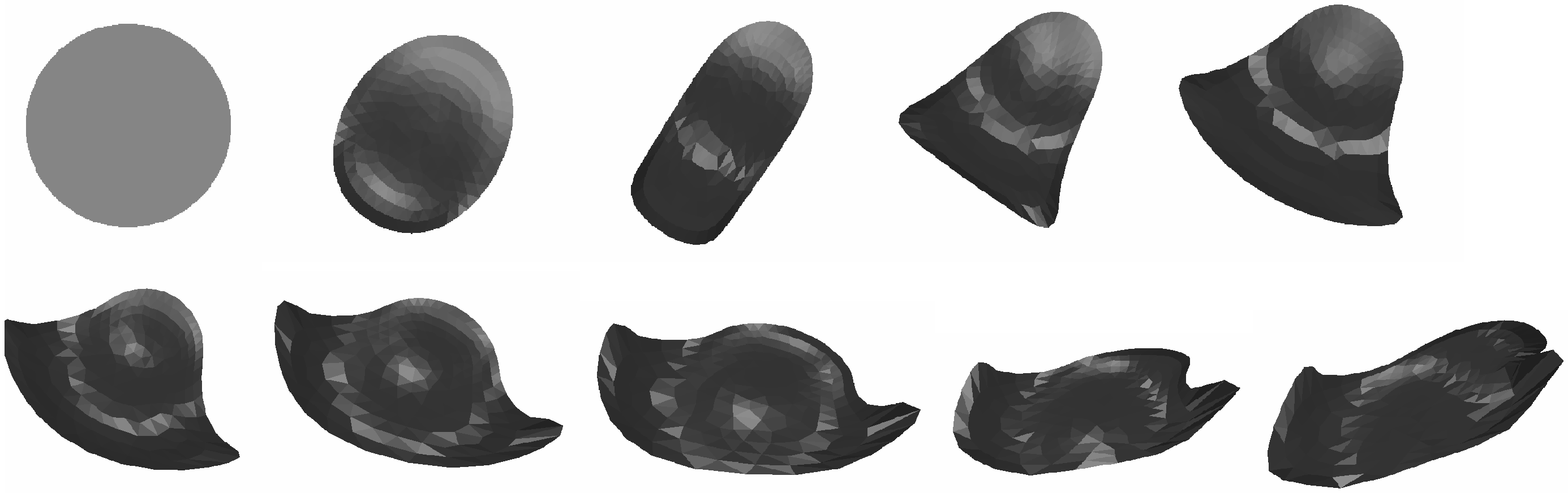}}
    \vskip -4truecm
    \caption{Evolution of the fluid domain}\label{Dom_evol1}
\end{figure}

We present some results of the simulation in figure
\ref{Dom_evol1}. Recalling that the shallow water equations are
based in the continuum mechanic, it cannot cut this domain with a
finite energy. Hence, even if the thickness of the layer of water
is very low, we conserve a thin film of water between all the
parts of the domain occupied by the fluid.
We recall that the computation is only made in the moving
two-dimensional domain $\Omega_t$. Mathematical study proves that
the main difficulty is to conserve the smoothness of the boundary.
This result can be observed in this simulation because, at the
final time, we have contact between two parts of the boundary (see
last part of the figure \ref{Dom_evol1}). Since we do not use a
``good'' boundary operator, we do not have sufficient smoothness
on the boundary. We observe then a change of connectedness and a
hole appears in the fluid domain.

\section{Conclusion}

The presented work suggests that capillarity effect needs to be
incorporated into the evolution equations, and more precisely in
the triple contact line. Theoretical results presented in
\cite{mateo} give sufficient smoothness for the capillarity term
used in equation \ref{tensionbord}.

The main difficulty is to describe the capillarity effects on the
depth integrated model and in a dynamical system. A possible
approach is to test some boundary operator and to compare
numerical and experimental result. This approach will be proposed
in future studies.

%



\end{document}